\documentclass[9pt]{article}
\usepackage{latexsym}
\usepackage{amssymb}
\usepackage[title]{appendix}
\usepackage{MnSymbol}
\usepackage{float}
\usepackage{booktabs}
\usepackage{mathrsfs}
\usepackage{amsmath}   
\usepackage{graphicx}
\usepackage{geometry}
\usepackage[font=footnotesize]{caption}

\begin{document}

\title {{\large \textbf{CHECKERBOARD GRAPH LINKS AND SIMPLY LACED DYNKIN DIAGRAMS}}}
\author{ \small{ LUCAS FERNANDEZ VILANOVA}}
\date{}
\maketitle
\renewcommand{\abstractname}{}
\begin{abstract}
{\footnotesize {\normalsize A}BSTRACT. We define an equivalence relation on graphs with signed edges, such that the associated adjacency matrices of two equivalent graphs are congruent over $\mathbb{Z}$. We show that signed graphs whose eigenvalues are larger than $-2$ are equivalent to one of the simply laced Dynkin diagrams: $A_{n}$, $D_{n}$, $E_{6}$, $E_{7}$ and $E_{8}$. Checkerboard graph links are a class of fibred strongly quasipositive links which include positive braid links. We use the previous result to prove that a checkerboard graph link with maximal signature is isotopic to one of the links realized by the simply laced Dynkin diagrams.\par}
\end{abstract}

\begin{center}
{\small 1. {\large I}NTRODUCTION} \\
\end{center}

Boileau, Boyer and Gordon showed that for strongly quasipositive links to have an L-space cyclic branched cover, they need to have maximal signature \cite{BBG1}, and more recently they showed that certain strongly quasipositive braids with a definite closure can be classified into the links realized by the simply laced Dynkin diagrams, abbreviated by $ADE$ diagrams. It is known that such classification holds for prime positive braid links with a positive Seifert form as proved by Baader, \cite{SB}. The notion of checkerboard graph links was introduced by S.\ Baader, L.\ Lewark and L.\ Liechti \cite{CGM} as a class of strongly quasipositive links strictly generalizing arborescent links with weights 2 and positive braids in a natural way. In \cite{CGM}, they ask if the $ADE$ diagrams correspond to the fibered links with maximal signature associated with checkerboard graphs. This is the essential motivation of this paper. Here is our main result: \\

\textbf{Theorem 1.1.}\textit{ A checkerboard graph link with maximal signature is isotopic to one of the links realized by the $ADE$ diagrams.}\\

For positive braid knots the topological $4$-genus is maximal exactly if the signature is maximal, \cite{genus}. Whether it is possible to obtain the same result for checkerboard graph links is also a question proposed by \cite{CGM}.\\

\textbf{1.2. Checkerboard graphs}\\

   A checkerboard graph is a finite, simple and plane oriented graph whose cycles are coherently oriented. The latter property is equivalent to say that they admit a checkerboard coloring, i.e., their dual, without the vertex corresponding to the \textit{unbounded} face, is a bipartite graph. The interest of studying these graphs resides in the fact that a checkerboard graph uniquely determines a strongly quasipositive fibred link [3, Theorem 2]. A special case of checkerboard graphs are the linking graphs, which uniquely determine a positive braid link [3, Theorem 1]. The manner to recover the link from a checkerboard graph or a linking graph is explained in detail in \cite{CGM}, see Figure 1 for an example of a positive braid link and the corresponding linking graph. In addition, as showed in \cite{CGM} it is possible to associate to a checkerboard graph an abstract open book i.e., a pair $(\Sigma, \phi)$, where $\Sigma$ is an oriented compact surface with boundary, and $\phi$ is a diffeomorphism, called the monodromy, that fixes the boundary pointwise. Baader and Lewark use the open books realized by checkerboard graphs in order to find two moves on these graphs that preserves the corresponding link type \cite{SL}, such moves will be of importance in the proof of Theorem 1.1.

A signed graph is a finite graph in which every edge is assigned a value $1$ or $-1$. Let $A(G)$ be its adjacency matrix (where the entries of the connected vertices are $+1$ or $-1$, depending on the sign of the connecting edge). Now, let $\beta$ be a positive braid word and $\Gamma (\beta)$ be its linking graph, then there exists a signed graph $\Gamma ^{\pm} (\beta)$ such that $2I+A(\Gamma ^{\pm} (\beta))$ is the symmetrized Seifert form of the closure of $\beta$, which we will denote as $L(\Gamma (\beta))$, (Proposition 1.4.2, \cite{Dehor}). 

\begin{figure}[H]
\includegraphics[scale=0.6]{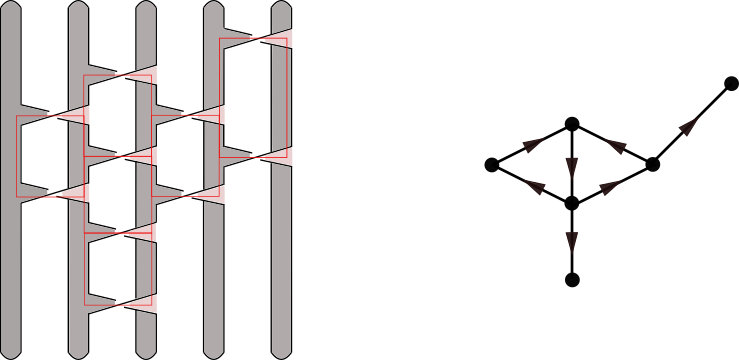}
\centering
\caption{ An example of the Seifert surface of a positive braid link and the corresponding linking graph. The red curves represent a natural homology basis. Note that they are in one-to-one with the vertices of the linking graph and two vertices are joined by an edge if the corresponding elements of the basis intersect.}
\end{figure}

This relates linking graphs and signed graphs, a similar relation can be made for checkerboard graphs  (such relation will be made more precise in Section 4). For that reason, we will start the proof of Theorem 1.1. by classifying positive symmetrized Seifert forms that we represent in terms of signed graphs.\\

\textbf{1.3. Signed graphs}\\

Let $G$ be a signed graph, we say that $G$ is a positive signed graph if $2I+A(G)$ is positive definite. For such graphs, we define a $t$-move that transforms one graph into another such that their adjacency matrices are congruent. This defines an equivalence relation that we call a $t$-equivalence.  \\

\textbf{Theorem 1.2.}\textit{ Let $G$ be a positive signed graph, then $G$ is $t$-equivalent to one of the $ADE$ diagrams.}\\

 A wide research in the field of spectral graph theory has been carried out on the graphs whose adjacency matrix have eigenvalues $>-2$. Results of the joint work of Cameron, Goethals, Seidel and Shult, \cite{CGSS} characterize those graphs \textit{represented} by one of the $ADE$ root systems. It is worth mentioning that a similar conclusion to the one in Theorem 1.2 can be achieved by using their results. However, the advantage of using the $t$-moves lies in the fact that, as we will further explain in Section 4, they are in close connection with certain checkerboard graph moves that preserve the link type.

If we only consider signed graphs that are planar and admit a checkerboard coloring (we shall call those graphs signed checkerboard graphs), we will find that Theorem 1.2. can be slightly sharpened. Indeed, we only need to use certain types of $t$-moves, which we call the $t'$-moves. We say that two graphs are $t'$-equivalent if there is a sequence of such moves relating one to another.\\

\textbf{Theorem 1.3.}\textit{ Let $G$ be a positive signed checkerboard graph, then $G$ is $t'$-equivalent to one of the $ADE$ diagrams.}\\

\begin{figure}[H]   
\includegraphics[scale=0.4]{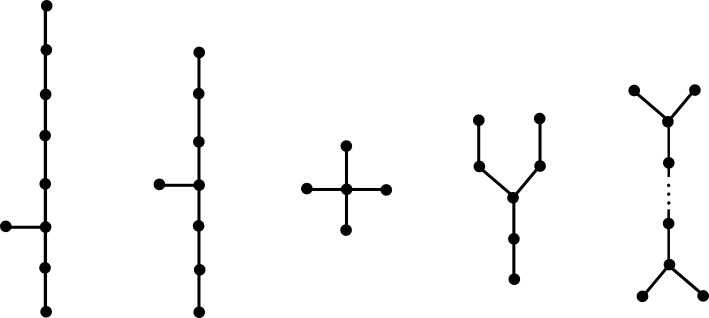}
\centering
\caption{ Forbidden minors, reading from left to right: $E$, $T$, $X$, $Y$ and $\tilde{D}$. Any signed graph containing an induced subgraph of these five types has an non-positive Seifert form.}
\label{fig:minors}
\end{figure}

 We will prove Theorems 1.2. and 1.3 by induction on the number of vertices. We show that if we add a vertex to one of the $ADE$ diagrams, we obtain a graph that is either $t'$-equivalent to one of the $ADE$ diagrams, or it has a non-positive Seifert form. To show the latter, we make use of the forbidden minors $E$, $T$, $X$ and $Y$, which do not have a positive definite symmetrized Seifert form, independently of their signs, \cite{SB}. Figure~\ref{fig:minors} shows the unsigned minors. In addition, we include the $\tilde{D}$ graph in our list of minors that have a positive semidefinite symmetrized Seifert form. A key step in the proof is that moves on unsigned graphs can be promoted to moves on signed graphs; therefore, simplifying the proof considerably.

\textbf{1.4. Outline}\\

 Section 2 provides a concise introduction to the $t$- and $t'$- moves on signed graphs and we give some of their properties. We will prove Theorems 1.2. and 1.3 in Section 3. In Section 4, we slightly generalize the moves proposed in \cite{SL} to find a version of the $t'$- moves for checkerboard graphs. We also explain the relation between sign graphs and checkerboard graphs, and we prove Theorem 1.1. \\

\textbf{1.5. Acknowledgements}\\

I thank Sebastian Baader for introducing me to the topics involved in this paper as well as his numerous and enlightening ideas concerning the proofs. I thank Lukas Lewark for explaining the checkerboard moves to me and also for the following discussions, which led to the proof of Lemma 4.1. I thank Sebastian Baader, Lukas Lewark and Livio Liechti for proofreading this paper. Also, I thank Arthur Bik for helpful comments.\\

\begin{center}
{\small 2. {\large M}OVES ON SIGNED GRAPHS} \\
\end{center}

We dedicate this section to define what we call a $t$-move on a signed graph, as well as the $t'$-moves, which are special cases of the former ones. Along that, we also study some of their properties that will be of importance in order to prove Theorems 1.2. and 1.3.

Before going into definitions, we discuss two three of detecting non-positive graphs aside from the forbidden minors mentioned above. Note that the latter are extremely useful when dealing with tree graphs, where signs can be ignored; however, once we encounter a cycle, the signs are important. Indeed, the following remarks show that the number of vertices and the number of negative edges in a cycle play a crucial role for detecting non-positive graphs. \\

\textbf{Remark 2.1.} Suppose $G$ is an $n$-cycle graph with $x$ number of negative edges and $n\geq 3$. Let $A(G)$ be its adjacency matrix. It is clear that if we change the signs of the two incident edges of a vertex in $G$, then the corresponding adjacency matrices are congruent. If the number of negative edges is even, we can transform the cycle into one with only positive edges, otherwise we can reduce the number of negative signs to one. Then, after a permutation of rows and columns 
\[
2I+A(G)\cong\left( \begin{array}{ccccccc}
2&1&& &  &\pm 1\\
1&2&1\\
 &1&2&1\\
 & & &\ddots\\
 & & &1&2&1\\
\pm 1 & & & &1&2\\
\end{array}\right).
\]
Where $\cong$ denotes matrix congruence over $\mathbb{Z}$. Here the entries $(1,n)$ and $(n,1)$ are positive if $x$ is even and negative otherwise. The determinant of the principal minor $(2I+A(G))_{n-a}$ is $n-a+1$ for $1 \leq a\leq n-1$, (see \cite{Meyer}). Hence, to determine whether $2I+A(G)$ is positive definite, it suffices to study its determinant. Using the cofactor expansion it is not hard to check that $det(2I+A(G))=0$ if $x$ and $n$ have the same parity and $det(2I+A(G))=4$ otherwise. \\

 When the number of negative edges and the number of vertices in a cycle have different parity we say that the cycle is positive. Remark 2.1 establishes a necessary condition for a signed graph with cycles to be positive.\\

\textbf{Remark 2.2.} Let $\Theta$ be a graph consisting of two positive cycles sharing $x\geq 2$ edges, $t$ of which are negative. Let $(n,p)$ and $(m,q)$ be these two cycles, where $n$ and $m$ are the number of vertices and $p$, $q$ stand for the number of negative signs in each cycle. Then, the \textit{outer cycle} is an induced subgraph with $m+n-2x$ edges and $p+q-2t$ of them are negative. Since the pairs $m,p$ and $n,q$ have different parity, it follows that the \textit{outer cycle} is not positive. Therefore, two positive cycles sharing more than one edge form a non-positive graph. \\

\begin{figure}[H]
\includegraphics[scale=0.35]{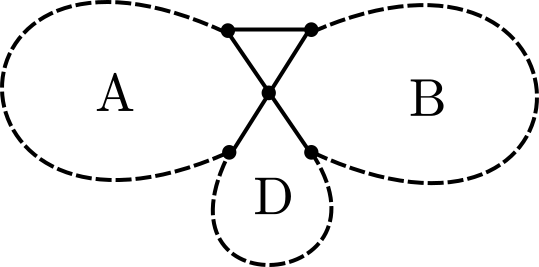}
\centering
\caption{ }
\label{fig:adb}
\end{figure}

\textbf{Remark 2.3.} Consider the graph in Figure \ref{fig:adb}, where $A$, $B$ and $D$ represent positive cycles of lengths $\geq 3$. The outer cycle, as an induced subgraph, is not a positive graph. The proof is similar to that in Remark 2.2. 

 We are now ready to define a $t$-move for signed graphs in three steps. Let $G$ be a signed graph, the cycles of which are positive:

\begin{itemize}
  \item[\textbf{Step 1.}] Pick any edge $\epsilon(x,y)$ in $G$, where $\epsilon\in\{ 1,-1\}$ is the sign of the edge, and choose one of its endpoints, say $x$.
  \item[\textbf{Step 2.}] Let $\{ v_{1},\dots,v_{n}\}$ be the set of vertices adjacent to $y$ (excluding $x$). In the case $x$ is the only adjacent vertex, jump directly to Step 3. Now, for all $v_{i}\in\{ v_{1},\dots,v_{n}\}$ draw an edge from $x$ to every vertex $v_{i}$, with the same sign as the edge $(y,v_{i})$ if $\epsilon = -1$ and opposite sign if $\epsilon = 1$. If an edge already exists, remove it. 
  \item[\textbf{Step 3.}] Change $\epsilon$ by $-\epsilon$.
\end{itemize}

We say that two signed graphs $G_{1}$ and $G_{2}$ are $t$-equivalent, and we denote it by $G_{1}\sim G_{2}$, if there exists a sequence of $t$-moves changing $G_{1}$ into $G_{2}$. Sometimes we write the pair $[v,w]$ to indicate that we perform a $t$-move on the edge $(v,w)$ and vertex $v$.

\begin{figure}[H]
\includegraphics[scale=0.6]{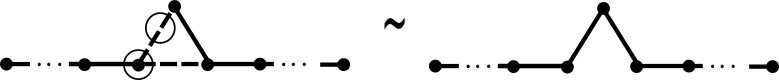}
\centering
\caption{An example of a $t$-move on the vertex and edge marked with a circle where the dashed lines indicate negative edges.}
\label{fig:1}
\end{figure}

Let $G$ be a signed graph and $|G|$ the graph we obtain from $G$ by ignoring the signs of its edges, this is usually called the underlying graph of $G$, a name that we adopt in this paper. Observe that the underlying graph of a signed graph that results from a $t$-move on $G$ does not depend on the signs of the edges in $G$. So it is possible to define the $t$-move for non-signed graphs (by simply ignoring edge signs in steps 1 and 2 and skipping step 3). \\

\textbf{Remark 2.4.} It is easy to verify that if we have a signed graph $G$ and a sequence of $t$-moves on $|G|$ such that $|G|\sim|G'|$, then the same sequence (choosing the same edges and vertices) transforms $G$ into $G'$, where the signs of $G'$ depend on those in $G$ and the chosen sequence. \\

Since we will consider forbidden minors that are trees and we  are interested on sequences that lead to tree graphs and the signs in a tree does not matter (see Remark 2.1); henceforth, we will consider underlying graphs only. This includes figures, starting at Figure \ref{fig:2}. 

For later use, consider the graph in Figure \ref{fig:2} (left side), which we call a B graph. It is not hard to check that:

\begin{figure}[H]
\includegraphics[scale=0.5]{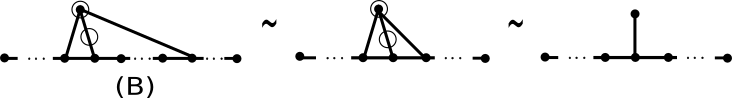}
\centering
\caption{ Graph B.}
\label{fig:2}
\end{figure}

Similarly, one can check the following relations:

\begin{figure}[H]
\includegraphics[scale=0.45]{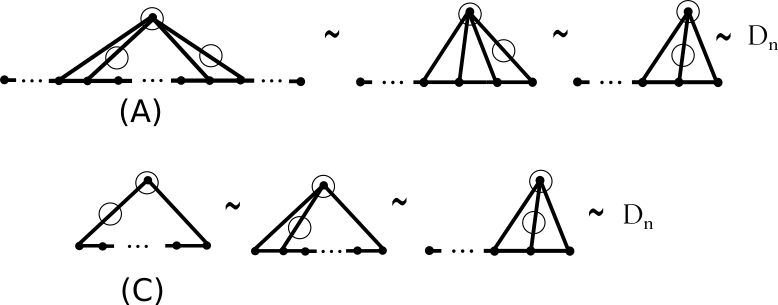}
\centering
\caption{ The graphs (A) and (C) are $t$-equivalent to $D_{n}$.}
\label{fig:3}
\end{figure}

The graphs $B$, $A$ and $C$ are not only an instructive example, but they will also be useful in the proof of Proposition 3.1.\\

\textbf{Lemma 2.3.}\textit{ Let $G_{1}$ be a signed graph the cycles of which are positive and let $G_{2}$ be a signed graph. If $G_{1}\sim G_{2}$, then $2I+A(G_{1})\cong 2I+A(G_{2})$.}\\

\textit{Proof:} Consider the vertex $v_{j}$ and the edge $(v_{i},v_{j})$ with sign $\epsilon$ in $G_{1}$. The matrix $2I+A(G_{1})$ has the form

\[
\left( \begin{array}{ccccccc}
   & &    & \vdots      &     & \vdots&   \\
   &  &  &  a_{mi}     &     & a_{mj}&   \\
     & &   &  \vdots     &     &  \vdots&  \\
    \dots & a_{im}  &\dots    &2 &   & \epsilon&\dots   \\
     &   &    &    &   &  & \\
     \dots& a_{jm}  &  \dots   &  \epsilon    &    &2& \dots \\
     &   &     &   \vdots   &   &\vdots&  \\
\end{array}\right).
\]

The row and column operations $R_{j}\rightarrow R_{j}\pm R_{i}$ and $ C_{j}\rightarrow C_{j}\pm C_{i}$, where we use the plus sign if $\epsilon =-1$ and the negative sign otherwise, give:

\[\left( \begin{array}{ccccccc}
   & &    & \vdots      &     & \vdots&   \\
   &  &  &  a_{mi}     &     & x_{mj}&   \\
     & &   &  \vdots     &     &  \vdots&  \\
    \dots & a_{im}  &\dots    &2 &   & -\epsilon&\dots   \\
     &   &    &    &   &  & \\
     \dots& x_{jm}  &  \dots   &  -\epsilon    &    &2& \dots \\
     &   &     &   \vdots   &   &\vdots&  \\
\end{array}\right).
\]

 Now, let $a_{mi}\neq 0$ and $a_{mj}\neq 0$, meaning that $v_{i}$ and $v_{j}$ are both connected to $v_{m}$ forming a $3$-cycle and since every cycle in $G$ is assumed to be positive, if $\epsilon=1$, then $a_{mi}=a_{mj}=\pm 1$ and if $\epsilon=-1$, then $a_{mi}=-a_{mj}$. Hence, after the row and column operation, $x_{mj}=0$ and $v_{j}$ loses its connection with $v_{m}$. If $a_{mi}\neq 0$ and $a_{mj}=0$, then $x_{mj}=a_{mi}$ for $\epsilon=-1$ and $x_{mj}=-a_{mi}$ for $\epsilon=1$, meaning that $v_{i}$ is now connected to $v_{m}$. If $a_{mi}=0$ and $a_{mj}\neq 0$, then $x_{mj}=a_{mj}$. Hence, the above matrix can be written as $2I+A(G')$ for some signed graph $G'$, where $G'$ is precisely the graph that results from performing the move on $G_{1}$ in the mentioned vertex and edge. Moreover, if $G_{1}$ has positive cycles, then the cycles of $G'$, if any, are also positive so if there is a sequence of moves changing $G_{1}$ into $G_{2}$, we can find a sequence of elementary operations changing $2I+A(G_{1})$ into $2I+A(G_{2})$. \hfill$\square$\\ 

It is clear now, that Theorem 1.2 and Lemma 2.3. implies that a positive definite matrix of the form $2I+A(G)$, for some signed graph $G$, is congruent over $\mathbb{Z}$ to a matrix $2I+A(\Gamma)$ where $\Gamma$ is one of the ADE diagrams.\\


\textbf{Definition 2.5.} A $t$-move on a vertex $v$ and edge $(v,w)$ where $deg_{w}\in\{1,2,3\}$ will be called a $t'$-move. We say that two signed graphs $G_{1}$ and $G_{2}$ are $t'$-equivalent if there exists a sequence of $t'$-moves changing $G_{1}$ into $G_{2}$. We denote it by $G_{1}\sim_{t'} G_{2}$. \\

For instance, all the moves in Figures \ref{fig:1} to \ref{fig:3} are $t'$-moves. For the case of positive, planar graphs that admit a checkerboard coloring the $t'$-moves can be restricted into the following three types (we have excluded the degree one case from Figure \ref{fig:4} for being trivial):

\begin{figure}[H]
\includegraphics[scale=0.45]{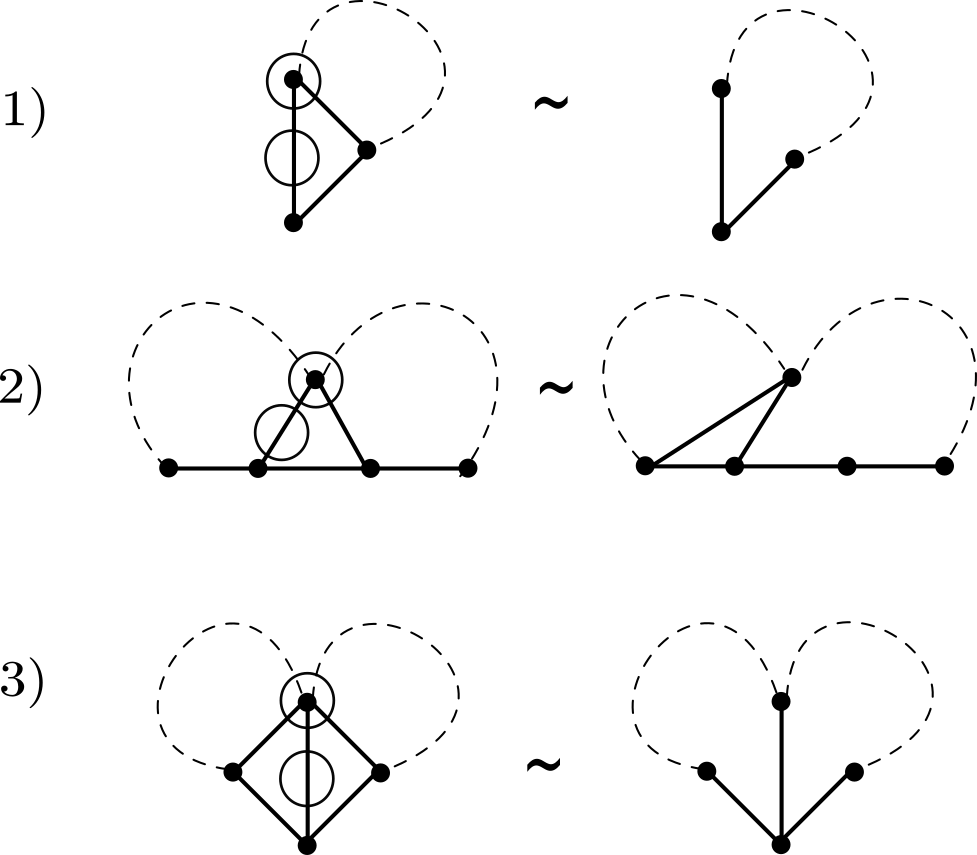}
\centering
\caption{ Dashed lines represents graphs connected to the line's endpoints.}
\label{fig:4}
\end{figure}

\begin{figure}[H]
\includegraphics[scale=0.65]{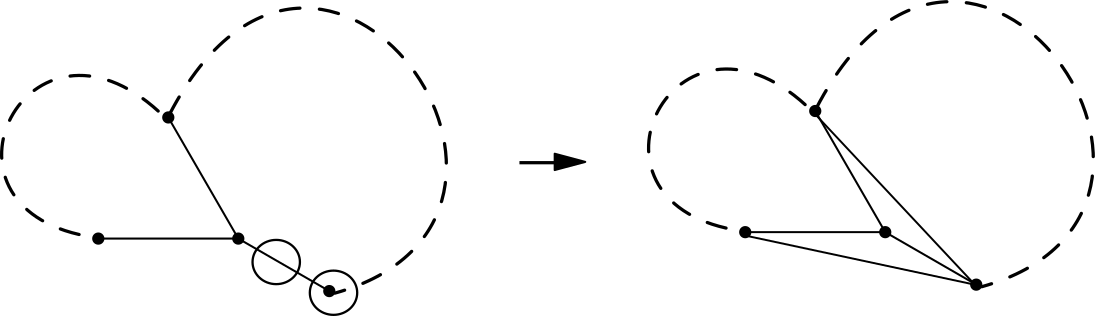}
\centering
\caption{}
\label{fig:relevant}
\end{figure}

To see that, note that positivity and checkerboard coloring properties on a planar signed graph $G$ implies that the maximum degree of a vertex, $v$, in $G$ is $6$. Moreover, $v$ is never an internal vertex (use Remark 2.3. and the fact that the wheel graph, $W_{7}$, is not positive). Also, exclude all the combinations that are $t'$-equivalent to a non-checkerborad graph, see for instance Figure \ref{fig:relevant} (such moves will not be allowed).

\begin{center}
{\small 3. {\large P}ROOF OF THEOREMS 1.2. AND 1.3}\\
\end{center}

In this section we first give the proof of Theorem 1.2. As explained in the introduction, the proof is done by induction. We assume that a graph, $\Gamma$, is $t$-equivalent to one of the $ADE$ diagrams, so it is clear that if we add a vertex $v$ to $\Gamma$, then $\Gamma \cup {v}\sim ADE \cup v$. Therefore, all we need to prove the Theorem 1.2 is to show that connecting a vertex to one of the $ADE$ diagrams results, after a sequence of $t'$-moves, into another $ADE$ diagram or a non-positive graph. We divide the proof in two parts: first, we study the case in which we add a vertex to the graph $A_{n}$ and second, we study in a similar manner those in which we add a vertex to $D_{n}$, $E_{6}$, $E_{7}$ and $E_{8}$. Finally, we prove Theorem 1.3. \\

\textbf{Proposition 3.1.}\textit{ If $G_{n}$ is a positive signed graph with $n$ vertices such that it is the union of $A_{n-1}$ with an extra vertex connected to $A_{n-1}$ by $m$ edges, then $G_{n}$ is $t'$-equivalent to one of the simply laced Dynkin diagrams.}\\

\begin{figure}[H]
\includegraphics[scale=0.6]{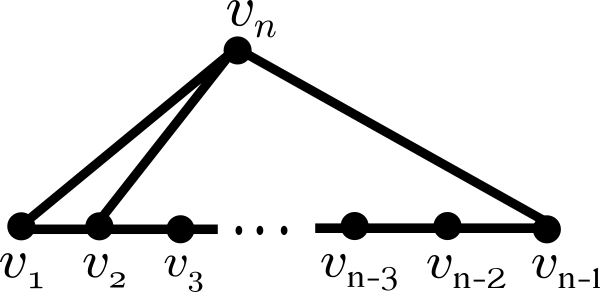}
\centering
\caption{ An example of $|G_{n}|$ for $m=3$.}
\label{fig:5}
\end{figure}

\textit{Proof:} The graph $|G_{n}|$ can be pictured as in Figure \ref{fig:5}, where $v_{n}$ has degree $m$ and the number of cycles in $G_{n}$ is therefore $m-1$.\\

Notice that, if $m>6$, then we can easily find an induced subgraph of type $X$ (see Figure \ref{fig:minors}) in $G_{n}$, implying that $G_{n}$ is not positive. Therefore, in order to prove the proposition, we consider all possible graphs for $m\leq 6$, which we divide in the following cases, and show that each of them is either $t'$-equivalent to a graph that contains a forbidden minor or it is $t'$-equivalent to one of the graphs $A_{n}$, $D_{n}$, $E_{6}$, $E_{7}$ or $E_{8}$.

First, observe that if $m=6$, then $G_{n}$ contains five cycles and these must have length $3$, otherwise $G_{n}$ contains a minor of type $X$. Hence, if we perform a $t'$-move on $v_{n}$ and the second and fifth edges we can reduce the degree of $v_{n}$ by four, which brings us to the $m=2$ case. In a similar fashion, If $m=5$, then there are four cycles in $G_{n}$. If at least two of these cycles have length $>3$, then $X\subset G_{n}$. So there must be at least three cycles of length $3$; one can easily check that after a $t'$-move on $v_{n}$ and one of the edges shared by these $3$-cycles, $deg_{v_{n}}=3$. Consequently,  we just need to consider the following four cases, in which $m$ takes the values $1$, $2$, $3$ and $4$.\\

\textbf{Remark 3.2:} Suppose that $v_{n}$ is connected to $v_{1}$ and $m>1$, i.e., $v_{n}$ is connected to at least one other vertex, say $v_{x}$. If $x=2$, then $[v_{n},v_{1}]$ reduces the degree of $v_{n}$ by one. If $x\neq 2$, then $[v_{n},v_{1}]$ creates a $3$-cycle and we can follow the sequence in Figure \ref{fig:2} in order to reduce the degree of $v_{n}$ by one, the same argument works in the case that $v_{n}$ is connected to $v_{n-1}$. Therefore, for the cases where $m>1$ we will assume that $v_{n}$ is not connected to any of these.\\

\textit{Case 1:} If $m=1$, and $v_{n}$ is connected to the vertex $v_{1}$ (or $v_{n-1}$), then $G_{n}\sim_{t'} A_{n}$. If $v_{n}$ is connected to the vertex $v_{2}$ (or $v_{n-2}$), then $G_{n}\sim_{t'} D_{n}$. If $v_{n}$ is connected to the vertex $v_{3}$ (or $v_{n-3}$), then for $n>8$; $E\subset G$, and for $n\in\{6,7,8\}$; $G\sim_{t'} E_{n}$. If $v_{n}$ is connected to a vertex different from those mentioned above and $n\geq 8$, then $T\subset G$.\\

\textit{Case 2:} If $m=2$, then there is one cycle in $G_{n}$, whose length we denote by $x$. If $x=3$, then $G_{n}\sim_{t'} A_{n}$, see Figure \ref{fig:1}. If $x=4$, then we obtain $D_{n}$ by a $t'$-move on $v_{n}$ and the two edges that have $v_{n}$ as an endpoint leading to a graph of type $B$ or $(a)$. For $x=5$ we need to consider the following case: 

\begin{figure}[H]
\includegraphics[scale=0.5]{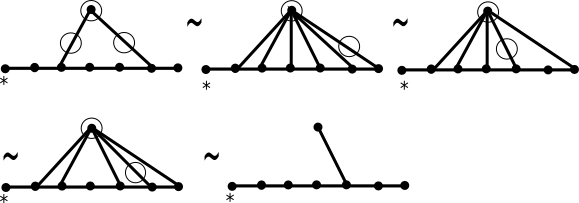}
\centering
\caption{ The non-signed graph without the starred vertex is $t'$-equivalent to $E_{7}$. The same $t'$-moves can be used with the additional starred vertex: leading to $E_{8}$.}
\end{figure}

\begin{figure}[H]
\includegraphics[scale=0.5]{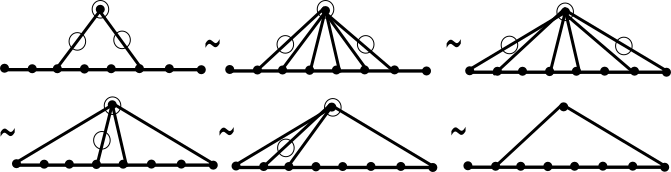}
\centering
\caption{ The non-signed graph is $t'$-equivalent to a graph that contains a $T$ minor. }
\label{fig:6}
\end{figure}

 From the graph in Figure \ref{fig:6}, it is clear that if we connect a new vertex to the starred one the resulting graph is $t'$-equivalent to $E$. If instead we connect a new vertex as it appears in Figure 8, then $G$ is $t'$-equivalent to a graph with a $T$ minor.

If $x>5$, then $G_{n}$ clearly contains a minor of type $\tilde {D}$.\\

\textit{Case 3:} For $m=3$, there are two cycles in $G_{n}$. If both of them have length $3$, we can perform a $t'$-move as in Figure \ref{fig:2} (middle case) which boils down to the $m=1$ case. If there is only one cycle of length $3$, then $G_{n}$ is a $B$ type graph and again, we are in case 1. If both cycles are $\geq 4$, then $\tilde{D}\subset G_{n}$.\\

\textit{Case 4:} Now, consider $m=4$, we know that there are three cycles in $G_{n}$ and if they all have length $>3$, then $X\subset G_{n}$. If there are two cycles with length $>3$ and the third has length $3$, then we can reduce $m$ by two as it is shown in Figure \ref{fig:7}. There are two cases to consider:
\begin{figure}[H]
\includegraphics[scale=0.5]{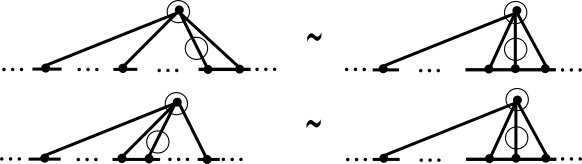}
\centering
\caption{ }
\label{fig:7}
\end{figure}
 If there are two non-adjacent cycles of length $3$ and one with length $>3$, then $G_{n}$ appears as in Figure \ref{fig:3}(A), so $G_{n}\sim_{t'} D_{n}$. If there are two adjacent $3$-cycles and one with length $>3$, then we can reduce $m$ by two, see Figure \ref{fig:7}, and if all three cycles have length $3$, then we can easily reduce $m$ by two, see the diagram of Figure \ref{fig:3}(A).\hfill$\square$\\

\textit{Proof of Theorem 1.2:} Let $n$ be the number of vertices of $G$. We may assume $n\geq 1$ and proceed by induction on $n$. For $n=1$ it is clear. Assume that the positive signed graph $G_{n-1}$ with $n-1$ number of vertices is $t'$-equivalent to one of the $ADE$ diagrams. In the following cases we show that if we connect one vertex to $G_{n-1}$, the new graph $G$ is either $t'$-equivalent to one of the $ADE$ diagrams or it is not positive (recall that we already know that this works for $G_{n-1}\sim_{t'} A_{n-1}$). Note that in a connected graph $G$ we can choose a vertex $v$ such that $G-v$ is connected (such a vertex exists; see e.g. \cite{Susanna}).\\

\begin{figure}[H]
\includegraphics[scale=0.55]{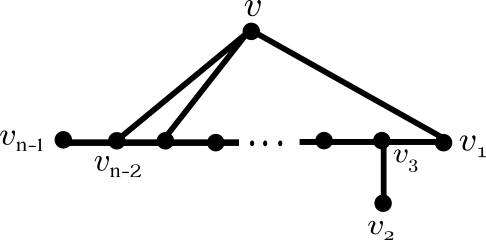}
\centering
\caption{ An example of $|G|$ for $m=2$.}
\label{fig:8}
\end{figure}

\textit{Case 1:} Let $G_{n-1}\sim D_{n-1}$ for $n>4$, so the graph $|G|$ can be pictured as in Figure \ref{fig:8}, where $v$ has degree $m$.

If $m>6$, then clearly $X\subset G$. Now recall that in the proof of Proposition 3.1 we only needed to consider the cases where $m\in\{1,2,3,4\}$, similar arguments work for this case. In addition, we observe that if $m>2$ and $v$ is connected to $v_{1}$ and $v_{2}$, then there are two cycles sharing two edges, which makes $G$ a non-positive graph for any choice of signs, see Remark 2.2. If that is not the case, then one can easily check that $m$ can be reduced to $1$, $2$ or $3$. Notice that, if $v$ is not connected to $v_{1}$, $v_{2}$ and $m>1$, then clearly $\tilde{D}\subset G$ or $X\subset G$ (except for $m=2$ and with a cycle of length $3$, in which case $G\sim_{t'} D_{n}$ by a $t'$-move). Thus, for $m>1$, we will only consider the cases when $v$ is connected to at least one of the vertices $v_{1}$ and $v_{2}$.  \\

\textit{Case 1.1:} If $m=1$ and $v$ is connected to the vertex $v_{1}$ or $v_{2}$ and $n\leq 8$, then $G\sim_{t'} E_{n}$ for $n\in\{6,7,8\}$. If $n>8$, then $G$ contains an induced subgraph of type $E$. If $v$ is connected to the vertex $v_{n-1}$, then $G_{n+1}\sim_{t'} D_{n+1}$. If $v$ is connected to any other vertices, then $\tilde{D}\subset G$ or $X\subset G$. \\

\textit{Case 1.2:} If $m=2$, the cycle has length 3 and $v$ is connected to $v_{1}$(or $v_{2}$) and $v_{3}$, then $G\sim_{t'} E_{n}$ for $n\in\{6,7,8\}$ by a $t'$-move on $[v,v_{1}]$ (or $[v,v_{2}]$), and it contains an $E$ minor for $n>8$. If the cycle has length $>3$ and $v$ is connected to $v_{i}$ and $v_{j}$ for $j=n-1$ and $i\in \{1,2\}$, then $G$ can be treated as one of the graphs in the proof of Proposition 3.1. But, if $3<j<n-1$ and $i\in \{1,2\}$, then one can check that:

\begin{figure}[H]
\includegraphics[scale=0.55]{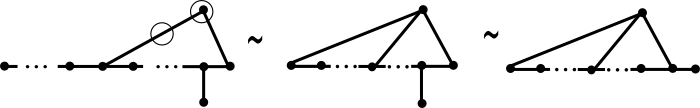}
\centering
\caption{}
\label{fig:9}
\end{figure}

 This again brings us to Proposition 3.1. Finally, if $v$ is connected to $v_{1}$ and $v_{2}$, then $G$ can also be treated as one of the graphs in the proof of Proposition 3.1.\\

\textit{Case 1.3:} If $m=3$, then $G$ has two cycles:

\begin{itemize}

 \item If both cycles have length $3$ and $v$ is connected to $v_{1}$, $v_{2}$ and $v_{3}$, then  

\begin{figure}[H]
\includegraphics[scale=0.6]{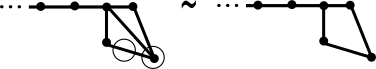}
\centering
\caption{ The non-signed graph is $t'$-equivalent to a graph as in the proof of Proposition 3.1, case 2.}
\end{figure}

 \item If there is only one cycle with length $3$, we need to consider either the case:

\begin{figure}[H]
\includegraphics[scale=0.5]{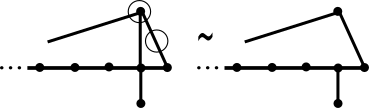}
\centering
\end{figure}

 or

\begin{figure}[H]
\includegraphics[scale=0.5]{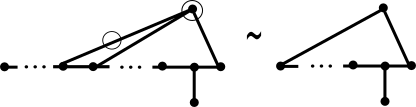}
\centering
\end{figure}

both of them reduce to case $1.2$. 

\item If both cycles have length $>3$, then we need to consider the following cases:

\begin{itemize} 

\item $v$ is connected to $v_{1}$, $v_{2}$ and $v_{i}$ for $3\leq i\leq n-1$, then $G$ has two cycles sharing two edges so by Remark 2.2, $G$ is not positive.
\item $v$ is connected to $v_{1}$, $v_{n-1}$ and $v_{i}$ for $3< i< n-1$, then $G$ is as in Figure \ref{fig:9} (middle case). 
\item If $v$ is connected to $v_{1}$ and two other vertices different from $v_{n-1}$, $v_{2}$ and $v_{3}$, then the only relevant case is when both cycles have length $4$ since for any other lengths one can easily find a minor of type $\tilde{D}$ or $X$. In the former case, Figure \ref{fig:10} shows that we can reduce it to case 1.2.

\begin{figure}[H]
\includegraphics[scale=0.5]{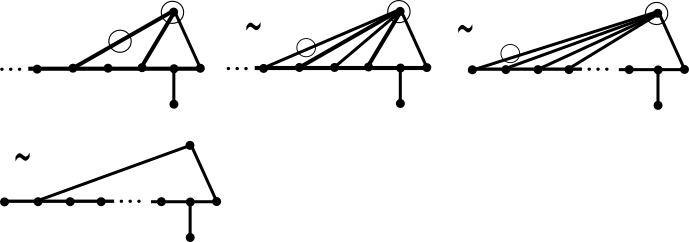}
\centering
\caption{ }
\label{fig:10}
\end{figure}

\end{itemize}
\end{itemize}

\begin{figure}[H]
\includegraphics[scale=0.55]{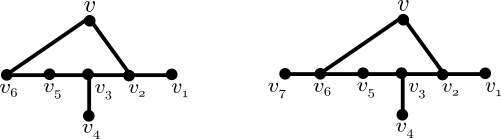}
\centering
\caption{ Examples of $|G|$ for $m=2$. Left: $G_{n-1}\sim_{t'} E_{6}$. Right: $G_{n-1}\sim_{t'} E_{7}$. We  can construct the $E_{8}$ case by simply connecting a new vertex, say $v_{8}$, to $v_{7}$. }
\label{fig:11}
\end{figure}

\textit{Case 2:} Assume $G_{n-1}$ is $t'$-equivalent to $E_{6}$, $E_{7}$ or $E_{8}$. So we can picture the graph $|G|$ as in Figure \ref{fig:11}.

\textit{Case 2.1:} First, consider that $G_{n-1}\sim_{t'} E_{6}$. If $m=1$ and $v$ is connected to the vertex $v_{1}$ or $v_{6}$, then $G\sim_{t'} E_{7}$. If $v$ is connected to the vertex $v_{3}$, then $X\subset G$. If $v$ is connected to a vertex different from those mentioned above, then $G$ contains an induced subgraph of type $\tilde{D}$. The case where $G_{n-1}\sim_{t'} E_{7}$ works similarly, and if $\sim_{t'} E_{8}$, then either we encounter the induced subgraphs $\tilde{D}$, $T$ and $X$, or $G=E$. \\

\begin{figure}[H]
\includegraphics[scale=0.35]{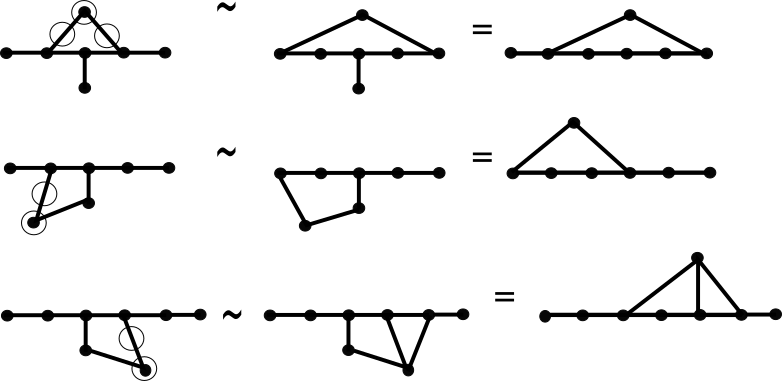}
\centering
\caption{}
\label{fig:13}
\end{figure}

\begin{figure}[H]
\includegraphics[scale=0.3]{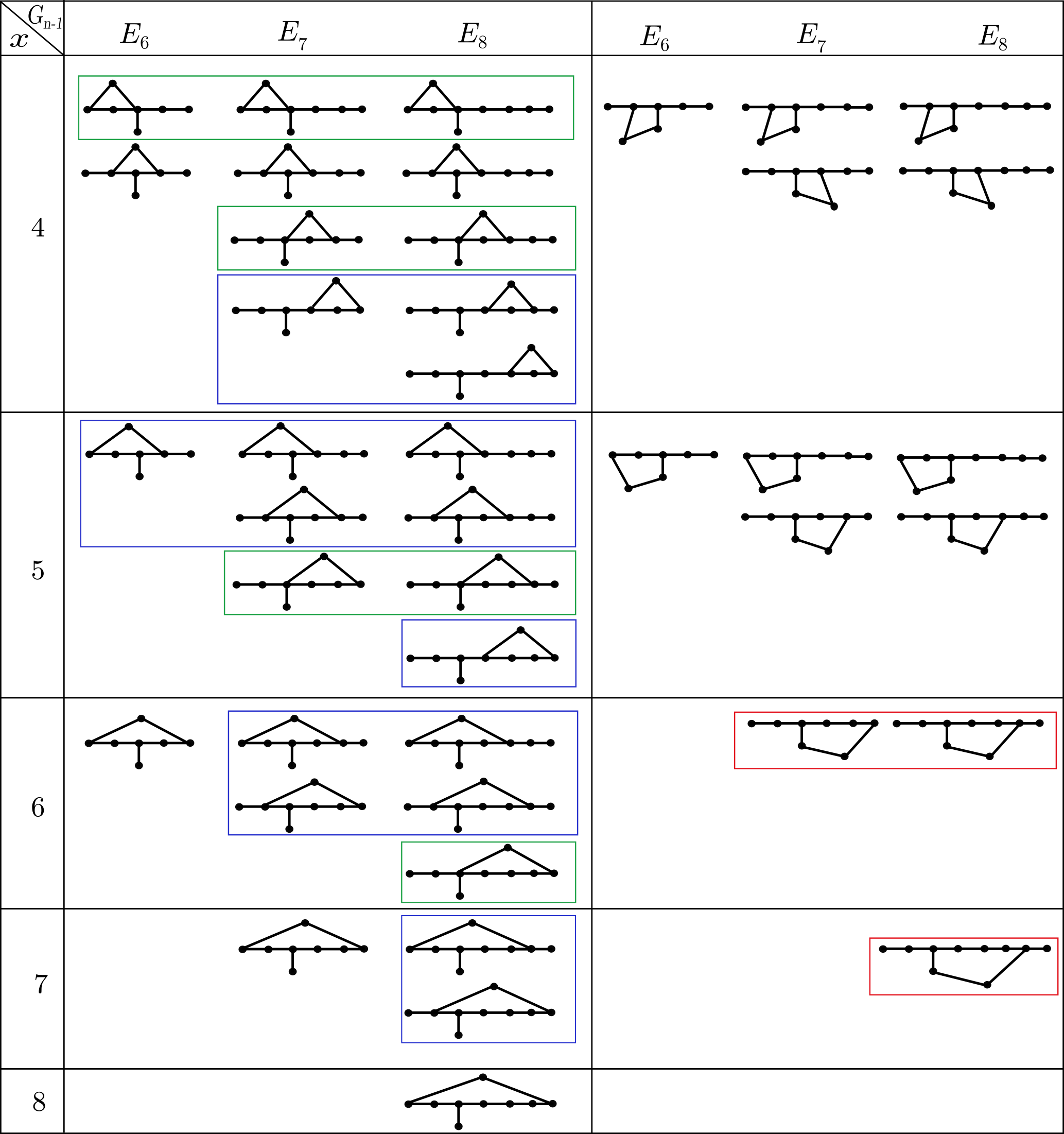}
\centering
\caption{ The graphs encircled with green contain an induced subgraph of type $X$, those encircled in blue contain $\tilde{D}$, and those in red contain $Y$. The right column includes the graphs in which $v$ is connected to $v_{4}$. The graphs that are not encircled are $t'$-equivalent to $E_{7}$, $E_{8}$ or, in the case where $G_{n-1}\sim_{t'} E_{8}$, to a graph with the induced subgraph $E$.}
\label{fig:12}
\end{figure}

\textit{Case 2.2:} Let $m=2$ so $G$ has one cycle of length $x$. First, if $x=3$ and $v$ is connected to $v_{3}$ and $v_{4}$, then $G\sim_{t'} Y$ (by a $t'$-move on $[v,v_{4}]$). If $G_{n-1}\sim_{t'} E_{6}$ and $v$ is connected to a different pair of vertices, then $G\sim_{t'} E_{7}$; however, that is not the case for $G_{n-1}\sim_{t'} E_{7}$, in which we can still find it is $t'$-equivalent to $E_{8}$ or it contains $T$. If $G_{n-1}\sim_{t'} E_{8}$, then it does not matter to which vertices we connect $v$, the resulting graph is $t$-equivalent to a graph with $T$ or $E$ as a minor. Now, for the cases where $x>3$ we simply draw all possible diagrams, see Figure \ref{fig:12}, and we find that they are either non-positive or $t'$-equivalent to $E_{7}$ or $E_{8}$ (see Figure \ref{fig:13}).

\textit{Case 2.3:} If $m=3$ and both cycles have length $3$, we can easily reduce the degree of $m$ by a $t'$-move, bringing us to the cases 2.1 or to the case 2.2 if $v$ is connected to $v_{3}$. If there is only one cycle of length $3$, then we can reduce the degree by one. If both cycles have length $>3$, one can check all possibilities as we did in Figure \ref{fig:12}, and find that all of them contain either $X$ or $\tilde{D}$.\\

Finally, if $m>3$ we can easily reduce the degree of $v$ by using $t'$-moves; taking us to the previous cases.\hfill$\square$\\

\textbf{Lemma 3.3.}\textit{ Let $G$ be a positive signed checkerboard graph with $n\geq 3$ vertices, then there exists at least one vertex in $G$ of degree $2$ or $3$.}\\

\begin{figure}[H]
\includegraphics[scale=0.35]{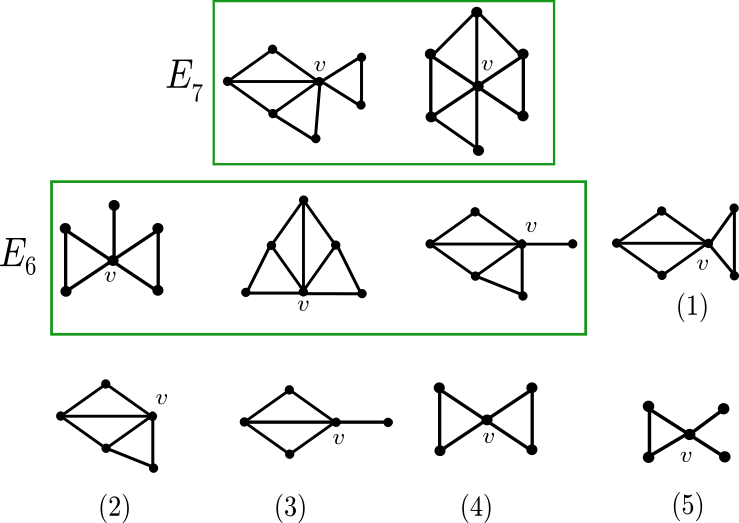}
\centering
\caption{}
\label{fig:4deg}
\end{figure}

\textit{Proof:} We show that if $G$ is a graph whose vertices have degree one or $\geq 4$, then $G$ is not positive. Recall from the previous section that positivity and checkerboard coloring properties on $G$ implies that the maximum degree of a vertex, $v$, in $G$ is $6$ and $v$ is never an internal vertex. Now, if $v$ is a vertex in $G$ of degree six, then the graphs on the first row in Figure \ref{fig:4deg} are the only two possible induced subgraphs involving $v$ in $G$ (any other combination is a non-positive graph). Similarly, if the degree of $v$ is five or four, we find that the graphs on the second and third row of Figure \ref{fig:4deg} are the possible induced subgraphs involving $v$.  

Theorem 1.2 shows that connecting a vertex to one of the special Dynkin graphs results into a graph in this class or into a non-positive one. Since the graphs in the green boxes of Figure \ref{fig:4deg} are $t'$-equivalent to $E_{6}$ or $E_{7}$, connecting more than two vertices to them will result into a non-positive graph. Thus, if we want to construct a positive graph whose vertices have degree one or $\geq 4$, we cannot use the graphs inside the boxes. 

\begin{figure}[H]
\includegraphics[scale=0.35]{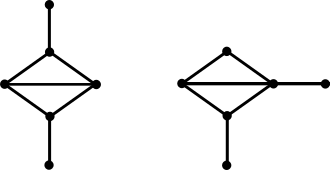}
\centering
\caption{ Two graphs $t'$-equivalent to $E_{6}$.}
\label{fig:e6}
\end{figure}

\begin{figure}[H]
\includegraphics[scale=0.35]{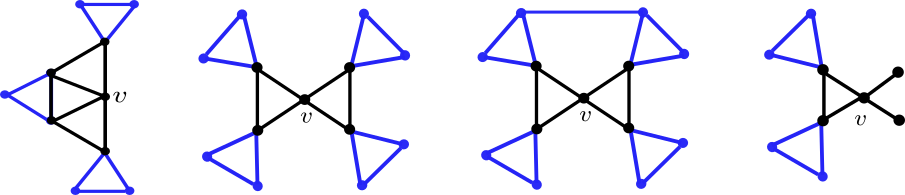}
\centering
\caption{}
\label{fig:cases}
\end{figure}

Now, consider the graphs in Figure \ref{fig:e6}, they are $t'$-equivalent to $E_{6}$. It is not hard to see that we cannot increase to $4$ the degree of all the adjacent vertices of $v$ in the graph (3), Figure \ref{fig:4deg} without encountering one of the graphs in Figure \ref{fig:e6} (neither in (1), since (3) $\subset$ (1)). Indeed, for the graphs (2), (4) and (5) only connecting triangles as it appears in Figure \ref{fig:cases} will work, but since the graph must have all its vertices of degree $\geq 4$ or one, for every vertex in the triangles we need to add at least two more vertices of degree one, creating a $\tilde{D}$ minor. Note that, joining two vertices of different triangles by an edge in the first graph of Figure \ref{fig:cases} results into a non-positive graph by Remark 2.3. By the same reason, we cannot join more than two triangles in the second graph. As for the fourth one, joining them results into a non-positive graph by direct computation.

\hfill$\square$\\

\textbf{Lemma 3.4.}\textit{ Let $\Gamma$ be a positive signed checkerboard graph. If there is a sequence of moves such that $\Gamma$ is $t'$-equivalent to one of the ADE diagrams, then $\Gamma \cup {v}$ is $t'$-equivalent to one of the ADE diagrams union $v$.}\\

\textit{Proof:} If there is a sequence of moves such that $\Gamma$ is $t'$-equivalent to one of the ADE diagrams and if $\Gamma \cup {v}$ is such that all the vertices in $\Gamma$ connected to $v$ have degree 2 or 3, then the same sequence can be used to transform $\Gamma \cup {v}$ into one of the ADE diagrams union $v$. If the degree  of a vertex connected to $v$ is $\geq 4$, then we cannot always use the same sequence. However, as we prove next, we can always find a $t'$-sequence that transform $\Gamma \cup {v}$ into one of the ADE diagrams union $v$.

Let $G$ be a positive planar graph such that it has a finite set of vertices of degree $\geq 4$, then there is a $t'$-sequence that transforms $G$ into a graph whose vertices have degree at most $3$. In order to prove it, assume that there exists in $G$ at least one vertex, say $w$, of degree $\geq 4$ that cannot be reduced by $t'$-moves to a degree $\leq 3$. Using Figure \ref{fig:4deg}, this means that at least three adjacent vertices of $w$ must have degree $\geq 4$ and cannot be reduced to lower degrees ($\leq 3$) either. By the proof of Lemma 3.3. we know that such graph is not positive. 

Now, let us come back to the case where $v$ is connected to a set of vertices of $\Gamma$, some of them with degree $\geq 4$. Then, we can find a sequence of $t'$-moves that transforms $\Gamma \cup {v}$ into $\Gamma' \cup {v}$ where $\Gamma'$ is positive and all its vertices have degree $\leq 3$. Thus, by Theorem 1.2. there is a sequence of $t$-moves transforming $\Gamma'$ into one of the ADE diagrams; but since all its vertices are of degree less or equal three, the statement follows. \hfill$\square$\\

\textit{Proof of Theorem 1.3:} We prove it by induction  on the number of vertices. Assume that there is a sequence such that $\Gamma$ is $t'$-equivalent to one of the ADE diagrams. By Lemma 3.4 we know that there is a $t'$-sequence transforming $\Gamma \cup {v}$ into one of the ADE diagrams union $v$, denote this graph by $G$. Now, analogous to the proof of Theorem 1.2 we get that $G$ is $t'$-equivalent to one of the ADE diagrams, completing the proof. \hfill$\square$\\

\begin{center}
{\small 4. {\large M}OVES ON CHECKERBOARD GRAPHS}\\
\end{center}

In this section we show that we can associate a signed graph to a checkerboard graph. We provide two checkerboard graph moves that preserve the corresponding link type; these are nothing but a generalization of the moves in [2]. Finally, we use these facts together with Theorem 1.3 to prove Theorem 1.1.

Recall from the definition of a checkerboard graph $\Gamma$, as it appears in \cite{CGM}, that this defines a strongly quasipositive fibred link, $L(\Gamma)$, constructed by plumbing positive Hopf bands according to the graph $\Gamma$. Thus, the Seifert matrix $V$ is an upper triangular matrix, \cite{einstein}. The signature, $\sigma (L(\Gamma))$, as defined by Trotter \cite{Trott}, is the signature of the symmetric matrix $M=V+V^{T}$. Thus, if $L(\Gamma)$ has maximal signature, $M$ must be a positive definite matrix. It follows that $M=[x_{ij}]$ is a symmetric matrix such that $x_{ij}=2$ if $i=j$. Now, since $x^{T}Mx>0$ for every non-zero vector $x\in \mathbb{R}^{n}$, choose $x=e_{i}\pm e_{j}$ then $(e_{i}+e_{j})^{T}M(e_{i}+e_{j})=2x_{ij}+4>0$ for all $1<i,j<n$ and $(e_{i}-e_{j})^{T}M(e_{i}-e_{j})=-2x_{ij}+4>0$, therefore

\[\mid x_{ij}\mid<2\ \ \ \ if\ \ \ i\neq j. \] 

and the non-diagonal coefficients are $0$, $1$ or $-1$. So $M$ is one of a finite list of matrices, and these can be represented by signed graphs. Moreover, since the vertices in a checkerboard graph and the ones in the corresponding signed graph represents positive Hopf bands, and two vertices are connected whenever the core curves of the corresponding Hopf bands intersect, it follows that the underlying graph of a checkerboard graph and the one of the corresponding signed graph are the same.\\

\textbf{Lemma 4.1.} \textit{The following two checkerboard graph moves preserve the corresponding link types.}

\begin{figure}[H]
\includegraphics[scale=0.7]{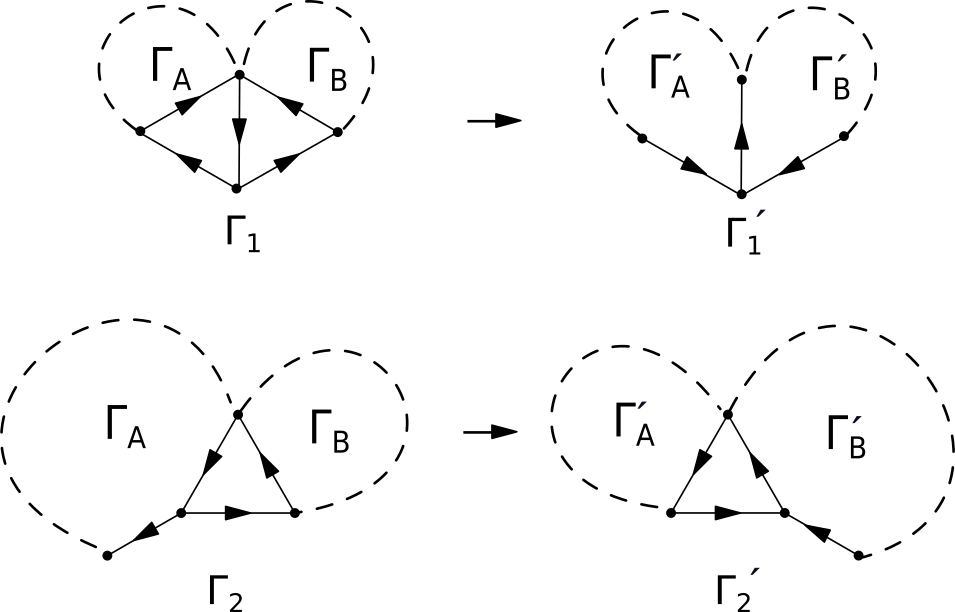}
\centering
\caption{$\Gamma _{A}$ and $\Gamma _{B}$ represents checkerboard graphs connected to the dashed line's endpoints.}
\label{fig:moves}
\end{figure}

\textit{Proof:} Recall that we can associate an abstract open book to a checkerboard graph. The goal is to show that the open books associated to the graphs $\Gamma _{1}$ and $\Gamma' _{1}$ are equivalent, i.e. let $(\Sigma _{1}, \phi _{1})$ and $(\Sigma' _{1}, \phi' _{1})$ be the open books associated to the graphs, then there is a diffeomorphism, $h$, between the surfaces $\Sigma _{1}$ and $\Sigma '_{1}$ such that $h\circ \phi '_{1}=\phi _{1}\circ h$. To do so, we will use the same argument as in \cite{SL}, and  check that these surfaces differ by a Dehn twist.

\begin{figure}[H]
 \includegraphics[scale=1.1]{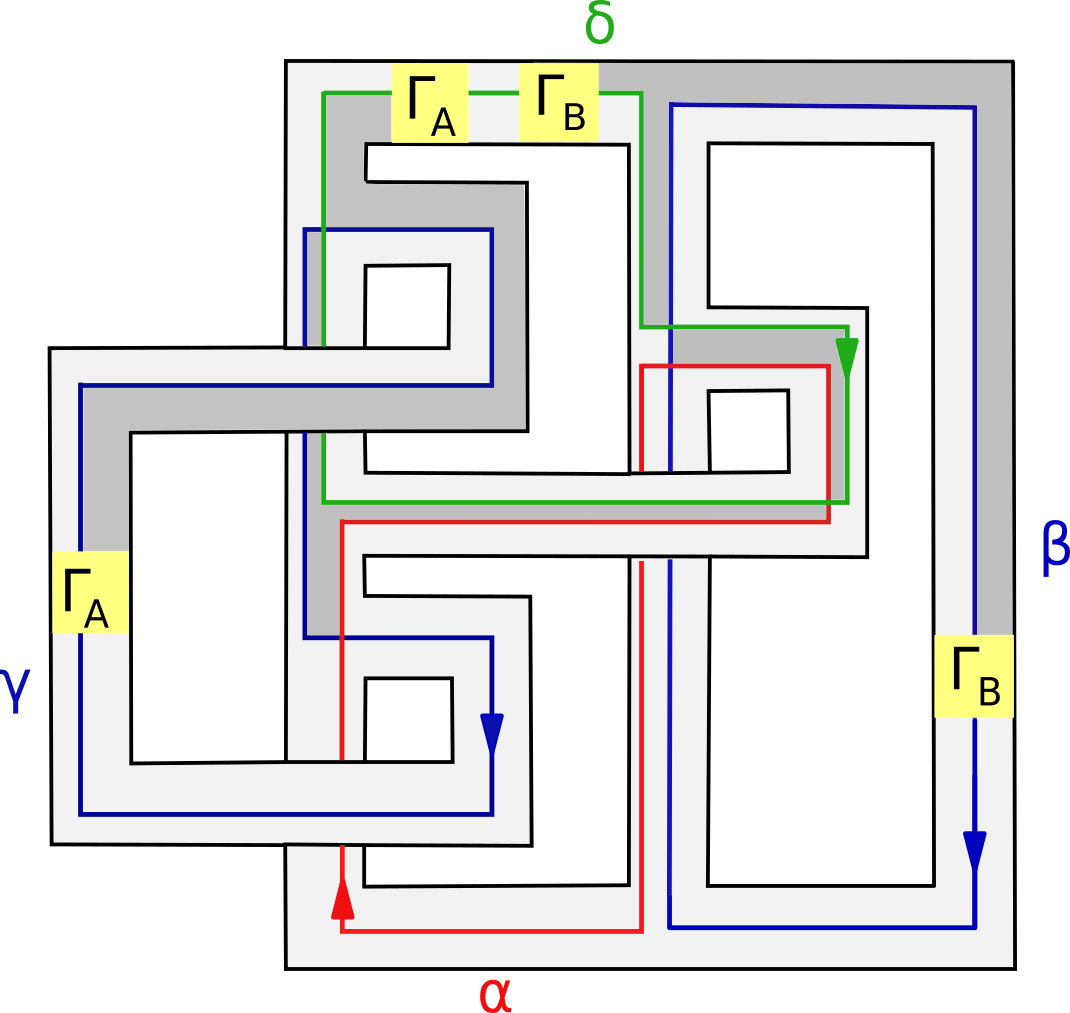}
 \includegraphics[scale=0.35]{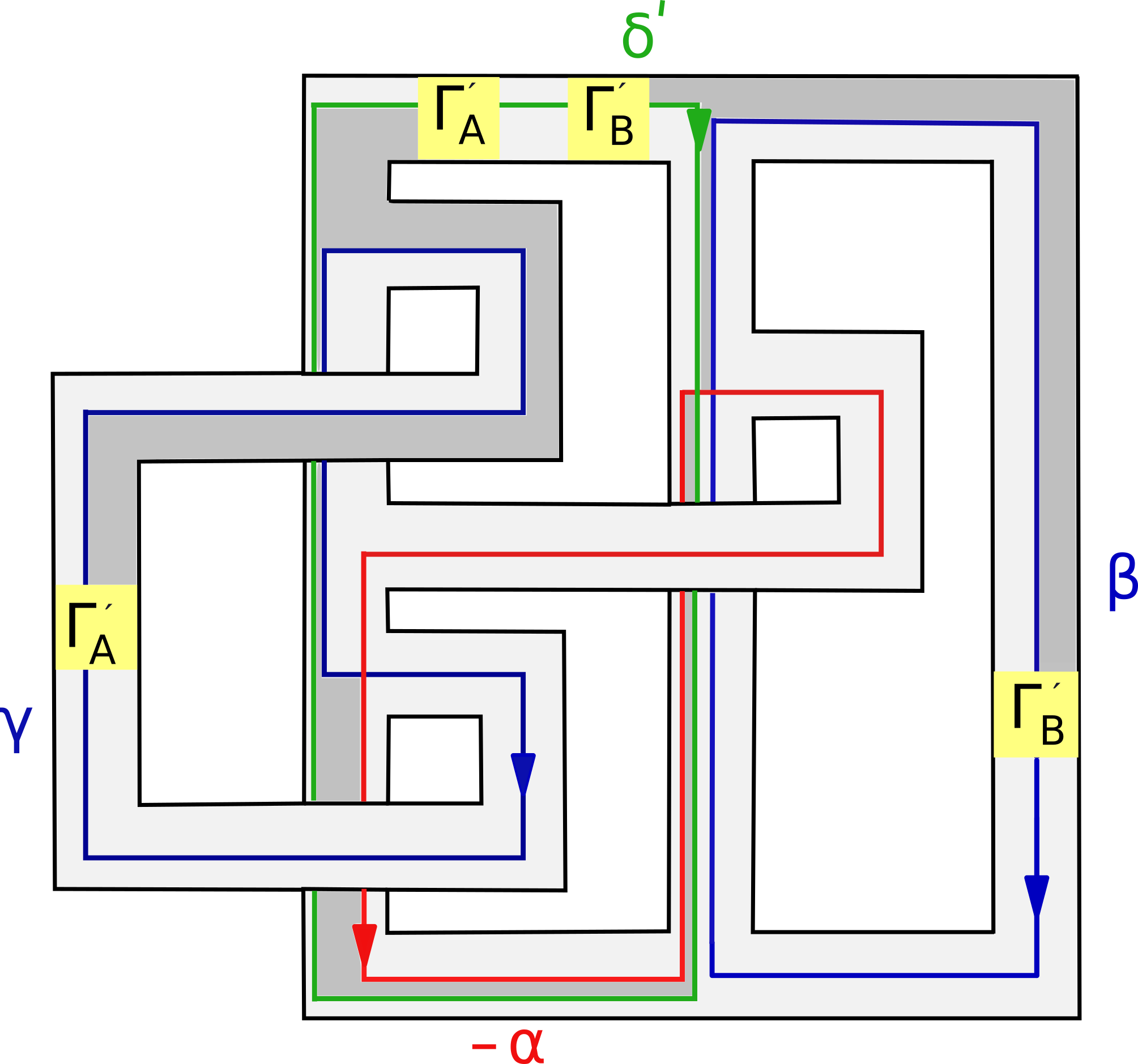}
 \centering
 \caption{ This drawing is a modification of a drawing in \cite{SL}. Right: four positive Hopf bands plumbed according to the graph in Figure \ref{fig:moves} with four $2$-handles(shaded regions). Left: the surface after a Dehn twist of $\delta$ along $\alpha$.}
 \label{fig:surf}
\end{figure}

The way we construct an abstract open book from a checkerboard graph is by gluing annuli; one for each vertex in $\Gamma$, and gluing disks; one for each cycle. The orientation of the core curves in each annulus is chosen so the intersection numbers with the other core curves corresponds to the orientations of the edges in $\Gamma$, see \cite{CGM} for further details about the construction. We can now construct the surface associated with $\Gamma _{1}$, see Figure \ref{fig:surf} (left), where the grey areas represents the disks and in the yellow squares we glue the parts that correspond to $\Gamma _{A}$ and $\Gamma _{B}$. Recall that the monodromy, $\phi$, is the product of positive Dehn twists in a certain order indicated by the orientation of the edges in $\Gamma _{1}$. In our case, if we label the vertices by $\alpha, \beta, \gamma, \delta$, one for each core curve in the surface and we take into account the orientation of the cycles $A$ and $B$, the monodromy can be written as $\phi _{A}\phi _{B}T_{\gamma}T_{\delta}T_{\alpha}T_{\beta}$ (it is customary to write $T_{a}$ as the Dehn twist of the curve $a$), where we have used the fact that we can switch two elements if there is no edge between them, \cite{CGM}. 

After performing a Dehn twist on $\delta$ along $\alpha$ we obtain the surface on the right, which one can easily check that corresponds to $\Gamma' _{1}$. Since $\delta'=T_{\alpha}^{-1}(\delta)$, then $T_{\delta'}=T_{\alpha}^{-1}T_{\delta}T_{\alpha}$ and $T_{\alpha}T_{\delta'}=T_{\delta}T_{\alpha}$. The monodromy is isotopic to $\phi _{A}\phi _{B}T_{\gamma}T_{\alpha}T_{\delta'}T_{\beta}$ which is precisely the monodromy that we obtain from the surface in the right. 

\begin{figure}[H]
\includegraphics[scale=1]{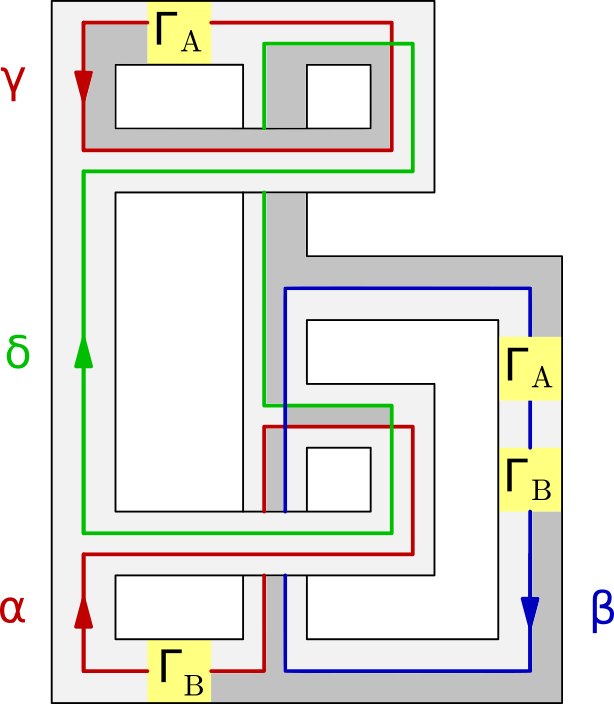}
\includegraphics[scale=1]{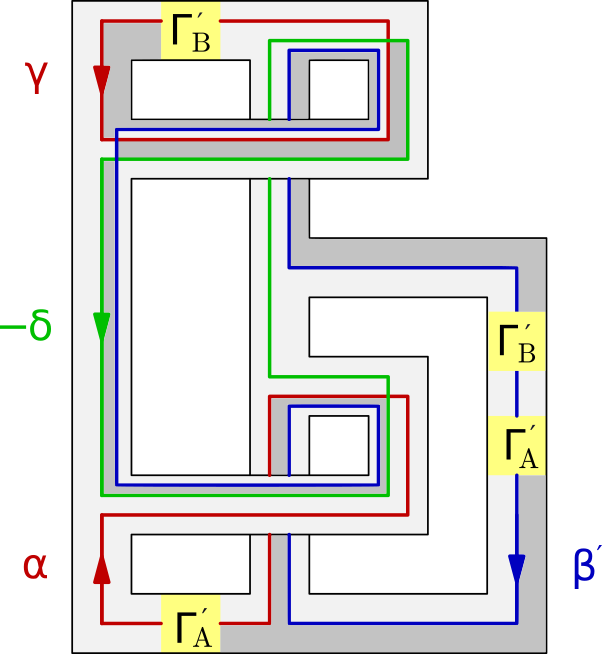}
\centering
\caption{Right: four positive Hopf bands plumbed according to the graph in Figure \ref{fig:moves} with three $2$-handles(shaded regions). Left: the surface after a Dehn twist of $\beta$ along $\delta$.}
\label{fig:surf2}
\end{figure}

 The second move is a generalization of the one described in \cite{SL}, only that this time we have the cycles $A$ and $B$, forming new grey regions as indicated in the corresponding abstract surfaces in Figure \ref{fig:surf2}. It is easy to check that the proof also works in this case.\hfill$\square$

\textit{Proof of Theorem 1.1:} Let $\Gamma$ be a checkerboard graph and let $L(\Gamma)$ be the associated link with maximal signature. We know that we can associate a signed graph, say $\Gamma^{\pm}$, to $\Gamma$ and by Theorem 1.3. there exists a sequence of $t'$-moves that transforms $\Gamma^{\pm}$ into one of the $ADE$ diagrams. Therefore, we can find a sequence of the moves in Lemma 4.1. (note that the moves in Figure \ref{fig:moves} are a checkerboard graph version of the $t'$-moves in Figure \ref{fig:4}) that transforms $\Gamma$ into a checkerboard graph of $ADE$ type preserving the corresponding link type. \hfill$\square$

Mathematics Institut, University of Bern, Sidlerstrasse 5, 3012 Bern, Switzerland\\
\textit{E-mail address:} lucas.fernandez@math.unibe.ch

\end{document}